\documentclass[12pt]{article} 
\usepackage{lipsum}
\usepackage[utf8]{inputenc}
\usepackage{hyperref}
\usepackage[T1]{fontenc}
\usepackage[a4paper, margin=2.7cm]{geometry}
\usepackage{amsmath, amsthm, amsfonts, amssymb}
\usepackage{mathrsfs}           
\usepackage{bbm} 
\usepackage{breakurl}
\usepackage{url}
\usepackage{authblk} 
\theoremstyle{definition}
\theoremstyle{remark}
\theoremstyle{example}
\usepackage{blindtext,titlefoot} 
\title{\bf On the notion of a quaternionic \\holomorphic function}
\author{Michael Parfenov}
\affil{\small \textit{Bashkortostan Branch of Russian Academy of Engineering, Ufa, Russia}} %DOBAVIL
 \date{}
\AtEndDocument{\bigskip{\footnotesize\noindent}}

\begin{document}
\maketitle\unmarkedfntext{\!\!\!\!\!\!\!\!\!\!\!2020 Mathematics Subject Classification: 30G35. Date: January 4, 2024  \\e-mail address: mwparfenov@gmail.com}
\begin{abstract}
A physically more adequate definition of a quaternionic holomorphic (H-holo\-morphic) function of one quaternionic variable compared to known ones and a quaternionic generalization of Cauchy-Riemann's equations are presented. \!\!At that a class of introduced H-holomorphic functions consists of those quaternionic functions whose left and right derivatives become equal after the transition to 3D space. \!\!The presented theory demonstrates a complete similarity of the algebraic properties and differentiation rules between the classes of  H-holomorphic and ordinary complex holomorphic functions, including the fact that quaternionic multiplication of the H-holomorphic functions behaves as commutative and the fact that each H-holomorphic function can be created from its complex holomorphic analogue by replacing a complex variable by a quaternion one. \!A fairly large number of detailed examples are given to illustrate the presented theory efficiency. \newline   
\end{abstract}
\section{Introduction}

The situation with regard to the derivative notion in quaternionic analysis seems somewhat confused. It has even led R. Penrose to say \!\cite[\!p.~\!201]{pen:road}, \textquotedblleft…quaternions are not really so mathematically \textquoteleft nice\textquoteright \,as they seem at first sigh. They are relatively poor \textquoteleft magicians\textquoteright; and, certainly, they are no match for complex numbers in this regard. The reason appears to be that there is no satisfactory quaternionic analogue of the notion of a holomorphic function.\textquotedblright 

During the last time there was a lot of interesting research (see, for example, \cite{sud:qua} - \cite{gss:rfqv}), including also new definitions of quaternionic holomorphy, however they do apparently not provide a full correspondence of properties of quaternionic holomorphic functions to complex analogues.

For a example, within the framework of existing conceptions, which consider most often the left or right derivative in isolation, regarding another as equivalent, we cannot construct a quaternionic holomorphic function from a corresponding complex holomorphic function of the same type by a direct replacement of a complex variable by a quaternion variable
in an expression for a complex holomorphic function although an analogous procedure is possible  (see, e.g., \cite[p.~353]{mf:mtp}) when constructing complex holomorphic functions from real-valued ones by a direct replacement of a real variable by a complex one.

\textit{In this regard, we’d like briefly to consider an essence of a simple approach based on physically more adequate notion of a quaternionic holomorphic function, which allows to call into question the above statement of R. Penrose and to establish a new class of quaternionic holomorphic functions with algebraic and differential properties fully similar to complex holomorphic analogues. Although this article has an overview character we would like to give a fairly large number of detailed examples illustrating the theory efficiency}.
\section{Quaternionic holomorphic functions}
The simple theory of quaternionic holomorphic functions can be built on principles fully similar (essentially adequate) to ones of complex holomorphic functions. At that the general definition of a derivative is based on the following main idea, viz.: each point of any real line is at the same time a point of some plane and 3D space as a whole, and therefore any characterization of differentiability (and its relations) at a point must be the same regardless of whether we think of that point as a point on the real axis or a point in the complex plane or a point in three-dimensional space. It follows that a quaternionic derivative of a quaternionic function\,
$\psi\left(p\right) $
must be defined similar to complex derivative as a limit of a difference quotient of $ \triangle\psi\left(p\right)$ by $\triangle p$ when
 $\triangle p$ converges to zero along any direction in the quaternionic space, where $\psi\left(p\right)=\psi_{1}\left(x,y,z,u\right)+\psi_{2}\left(x,y,z,u\right)i+\psi_{3}\left(x,y,z,u\right)j+\psi_{4}\left(x,y,z,u\right)k
$ is a quaternionic function of real components $x,y,z,u$ of an independent quaternionic variable $p=x+yi+zj+uk$.
The functions $\psi_{1}\left(x,y,z,u\right),\psi_{2}\left(x,y,z,u\right),\psi_{3}\left(x,y,z,u\right),\psi_{4}\left(x,y,z,u\right)$ are real-valued and $i,j,k$ are the base quaternions of the quaternion space $\mathbb{H}$.

Using the Cayley–Dickson construction (doubling form) \cite{ks:hn} we have
\begin{equation} \label{Eq1} 
 \begin{gathered}
p=a+b\cdot j\in \mathbb{H},\;\;\;\;\;\;\psi\left(p\right)=\Phi_{1}\left(a,b,\overline{a}, \overline{b} \right)+\Phi_{2\cdot}\left(a,b,\overline{a}, \overline{b} \right)\cdot j \in \mathbb{H},  \\
\,\,\, \overline{ p}=\overline{ a}-b\cdot j\subset \mathbb{H}\text{,} \;\;\;\;\;\; \overline{ \psi\left(p\right)}=\overline{\Phi_{1}\left(a,b,\overline{a}, \overline{b} \right)}\;- \Phi_{2\cdot}\left(a,b,\overline{a}, \overline{b} \right)\cdot j \subset \mathbb{H} \text{,}
\end{gathered}
\end{equation}
where
\begin{equation}  \label{Eq2}
 \begin{gathered}
 a=x+yi,\;\;\;\;\;b=z+ui, \\ 
 \overline{a}=x-yi,\;\;\;\;\; \overline{b}=z-ui,
\end{gathered}
\end{equation}
$$\Phi_{1}\left(a,b,\overline{a}, \overline{b} \right)=\psi_{1}+\psi_{2}i, \;\;\;\;\;\Phi_{2}\left(a,b,\overline{a}, \overline{b} \right)=\psi_{3}+\psi_{4}i,$$
$$\overline{ \Phi_{1}\left(a,b,\overline{a}, \overline{b} \right)}=\psi_{1}-\psi_{2}i,\;\;\;\; \overline{ \Phi_{2}\left(a,b,\overline{a}, \overline{b} \right)}=\psi_{3}-\psi_{4}i$$
are compex quantities, the \textquotedblleft  $\cdot$\textquotedblright  \,and overbar signs denote, respectively, quaternionic multiplication and complex (or quaternionic if needed) conjugation. For simplicity, we use also the short designations $\Phi_{1}$ and $\Phi_{2}$ respectively instead of $\Phi_{1}\left(a,b,\overline{a}, \overline{b} \right)$ and $\Phi_{2}\left(a,b,\overline{a}, \overline{b} \right)$.

Since the quaternion algebra is a noncommutative algebra with division \cite{ks:hn}, there can exist the left and the right definition of a quaternionic derivative:
$$\psi_{left}'=\lim_{\triangle p \rightarrow 0}\left[\left(\triangle p\right)^{-1}\cdot\left\{ \psi\left(p\,+\triangle p \right)-\psi\left(p\right)\right\}\right]\!,
$$
$$\,\,\,\psi_{right}'=\lim_{\triangle p \rightarrow 0}\left[\left\{\psi\left(p\,+\triangle p \right)-\psi\left(p\right)\right\}\cdot\left(\triangle p\right)^{-1}\right]\!.
$$
Usually, one of them is used (regarding another as equivalent \cite{sud:qua,leo:quat}, etc). However such an approach cannot be recognized as full, since the algebra underlying the quaternionic analysis  needs at the same time to use both quaternionic multiplications: on the left and on the right (see details in\! \cite[\!\!p.\!\!~32]{ks:hn} and\!  \cite[\!\!p.\!\!~7]{pm:aqg}) to represent \textit{all arbitrary (and noncommutative) rotations of any vector in 3D space} just like as the complex algebra  underlying the complex analysis uses the complex  multiplication, representing all  arbitrary (commutative) rotations of any vector in the complex plane.  Thus  we \textit{have to} refuse to consider only the left or only the right approach (regarding another as equivalent) when defining a quaternionic derivative. The left and right approaches should be considered only together \cite{pm:aqg}.

In complex analysis the notion of a derivative of a holomorphic function (as a complex potential) is always associated \cite{mf:mtp,mh:ca}  with \textit{an unambiguous strength value} of some physical planar vector field (electrical field, fluid flow and etc.) and accordingly \textit{the derivative is unambiguous}. Similarly, the quaternionic derivative must correspond some physical field of dimension 3 and hence  \textit{must be also unambiguous. }Thus we \textit{are forced to require} an equality of the left and right derivatives, i.e.$\;\psi_{left}'=\psi_{right}'$ in a domain of a function definition  to associate the quaternionic derivative with physical reality of 3D space. 

In fact, we want to define the notion of a quaternionic differentiability, which will be \textquotedblleft physically more adequate \!\!\textquotedblright \,(or essentially adequate) than the left or the right approach in isolation. At first sight the requirement $\psi_{left}'=\psi_{right}'$\;taken alone seems absurd, but we will try to establish those quaternionic functions and those conditions for which it is fair. As we see below, such the functions are desired quaternionic holomorphic functions and such the conditions are quaternionic Cauchy-Riemann's equations that are performed for these functions  after the transition to 3D space. At that  the left and the right derivatives of such functions become equal  after the transition to 3D space.

According to these principles, it follows that the limit of the difference quotient is required to be independent not only of directions to approach a limiting point when $\triangle p \rightarrow 0$ (as in complex analysis \cite{mh:ca,mf:mtp}), but also of the manner of quaternionic division: on the left or on the right.  In this sence we can say that  \textquotedblleft the derivative is independent of the way of its computation \textquotedblright. Based on that, it is possible to introduce the following\\\\ \textbf{Definition 1.}  \textit{A single-valued quaternionic function $\psi\left(p\right):G_{4}\rightarrow \mathbb{H}$ is quaternionic differentiable at a point $p\in G_{4}\subseteq \mathbb{H}$ if there exists a limiting value of the difference quotient ${\triangle \psi\left(p\right)}/{\triangle p}$ as $\triangle p\rightarrow 0,$ and this value is independent of the way of its computation.} \\\\As we shall see further, such an independence is provided after the transition to 3D space.\\\\
 \textbf{Definition 2.} \textit{A quaternionic function is said to be a quaternionic holomorphic function (briefly, H-holomorphic function) at a point $p\in\mathbb{H}$  if it has a quaternionic derivative independent of a way of its computation in some open connected neighborhood $G_{4}\subset \mathbb{H}$ of a point $p$}.
 \\\\In the Cayley–Dickson doubling form Definitions 1 and 2 (by equating the left and the right derivatives) lead to the following formulation of the necessary and sufficient conditions \cite[\!p.p.\!~14-19, p.~22]{pm:aqg} for the function $\psi\left(p\right)$ to be H-holomorphic: \\\\
 \textbf{Definition 3.} \textit{It is assumed that the constituents $\Phi_{1}\left(a,b,\overline{a}, \overline{b} \right)$ and $\Phi_{2}\left(a,b,\overline{a}, \overline{b} \right)$ of a quaternionic function $\psi\left(p\right)=\psi\left(a,b,\overline{a}, \overline{b} \right)=\Phi_{1}+\Phi_{2}j$ possess continuous first-order partial derivatives with respect to $a,b,\overline{a}$ and $\overline{b}$ in some open connected neighborhood  $G_{4}\subset \mathbb{H}$ of a point $p\in\mathbb{H}.$ Then a function $\psi\left(p\right)$ is said to be H-holomorphic and denoted by $\psi_{H}\left(p\right)$ at a point $p,$ if and only if the functions  $\Phi_{1}\left(a,b,\overline{a}, \overline{b} \right)$ and $\Phi_{2}\left(a,b,\overline{a}, \overline{b} \right)$ satisfy in $G_{4}$ the following quaternionic generalization of complex Cauchy-Riemann's equations} \cite{pm:aqg}:
\begin{equation} \label{Eq3} \left\{
\begin{aligned}
1)\,\,\,\,(\,\partial_{a}\Phi_{1}\!\!\mid \,\,=(\,\partial_{\overline{b}}\overline{\Phi_{2} }\! \mid, \,\,\,\,\,\,\,\,\,\,\,\,\,\,2)\,\,\,\,(\,\partial_{a}\Phi_{2}\!\!\mid &=-\,(\,\partial_{\overline{b}}\overline{\Phi_{1} }\! \mid 
,\\
3)\,\,\,\,(\,\partial_{a}\Phi_{1}\!\!\mid \,\,=(\,\partial_{b}\Phi_{2}\! \mid, \,\,\,\,\,\,\,\,\,\,\,\,\,\,4)\,\,\,\,(\,\partial_{\overline{a}}\Phi_{2}\!\!\mid &=-\,(\,\partial_{\overline{b}}\Phi_{1}\!\! \mid .
\end{aligned}
\right. 
\end{equation}
Here $\partial_{i},\,i=a,\overline{a},b,\overline{b}$ denotes the partial derivative with respect to $i$. The brackets $\left(\cdots{}\!\!\mid\right.$ with the closing vertical bar indicate that the transition $a=\overline{a}=x$ (to 3D space) has been already performed in expressions enclosed in brackets. Equations (\ref{Eq3}-1) and (\ref{Eq3}-2) have the components relating to the left quaternionic derivative and equations (\ref{Eq3}-3) and (\ref{Eq3}-4) to the right one \cite[p.p.~16,17]{pm:aqg}. The requirement  $a=\overline{a}=x$ provides the possibility of joint implementation of equations (\ref{Eq3}-1,2) for the left quaternionic derivative and (\ref{Eq3}-3,4) for the right one. Equations (\ref{Eq3}-1) and (\ref{Eq3}-3) as well as (\ref{Eq3}-2) and (\ref{Eq3}-4) become, respectively, identical after the transition to 3D space \!\cite[p.~18]{pm:aqg}, i.e. the left derivative becomes equal to the right one.

We see that H-holomorphy conditions (\ref{Eq3}) are defined so that during the check of the quaternionic holomorphy of any quaternionic function we have to do the transition  $a=\overline{a}=x$ in already calculated expressions for the partial derivatives of the functions $\Phi_{1}\left(a,b,\overline{a}, \overline{b} \right)$ and  $\Phi_{2}\left(a,b,\overline{a}, \overline{b} \right)$ and their complex conjugations. However, this doesn’t mean that we deal with triplets in general, since the transition $a=\overline{a}=x$ (or $y=0$) cannot be initially done for quaternionic variables and functions. Any quaternionic function remains the same 4-dimensional quaternionic function regardless of whether we check its holomorphy or not. Simply put, the H-holomorphic functions are 4-dimensional quaternionic functions for which the partial derivatives of components of  the Cayley–Dickson doubling form satisfy equations (\ref{Eq3}) after the transition to 3D space. In other words, they are those quaternionic functions whose the left and the right derivatives become equal after the transition to 3D space. That's not surprising, since unambiguous stationary physical fields, represented by derivatives, exist precisely in 3D space. 
\subsection{Example of \mathversion{bold}$p$ to the power of 2}  \label{Exa1}
  \textit{The function $\psi\left(p\right)=p^{2}=\left(a+b\cdot j\right)\cdot \left(a+b\cdot j\right)=\Phi_{1}\left(a,b,\overline{a}, \overline{b} \right)+\Phi_{2}\left(a,b,\overline{a}, \overline{b} \right)\cdot j. \, \,$}By direct quaternionic multiplication we obtain the following expressions for the components of the Cayley–Dickson doubling form: $$\Phi_{1}\!\!\left(a,b,\overline{a}, \overline{b} \right)=a^{2}-b\overline{b}\,\,\,\,\,\,\,\, \text{and}\,\,\,\,\,\,\,\,\Phi_{2}\!\!\left(a,b,\overline{a}, \overline{b} \right)=\left(a+\overline{a}\right)b.$$ Correspondingly, the complex conjugate functions are $$\,\,\overline{ \Phi_{1}\left(a,b,\overline{a}, \overline{b} \right)}=\overline{a}^{2}-\overline{b}b\,\,\,\,\,\,\,\, \text{and}\,\,\,\,\,\,\,\,\overline{ \Phi_{2}\left(a,b,\overline{a}, \overline{b} \right)}=\left(\overline{a}+a\right)\overline{b}.$$ When extracting the components $\Phi_{1}\!\!\left(a,b,\overline{a}, \overline{b} \right)$ and $\Phi_{2}\!\!\left(a,b,\overline{a}, \overline{b} \right)$ we always take into consideration the fact that $j\alpha=\overline{\alpha}j$ \cite[p.~42]{ks:hn} for any $\alpha\in\mathbb{C}$, where $\mathbb{C}$ is the complex plane and $j$ is the imaginary quaternion unit, so its square is -1.
 
 Now we calculate the derivatives: $\partial_{a}\Phi_{1}=2a,\,\,\,\partial_{\overline{b}}\overline{\Phi_{2} }=\overline{a}+a,\,\,\,\partial_{a}\Phi_{2}=~b,\,\,\,\partial_{\overline{b}}\overline{\Phi_{1} }=-b$ for equations (\ref{Eq3}-1,2) and $\,\,\,\partial_{b}\Phi_{2} =a+\overline{a},\,\,\,\partial_{\overline{a}}\Phi_{2}=b,\,\,\,\partial_{\overline{b}}\Phi_{1}=-b$ for equations (\ref{Eq3}-3,4). By substituting $a=\overline{a}$, i.e. \,$a=x$ and \,$\overline{a}=x$ into these
expressions for partial derivatives we have
\begin{equation*}  \left\{
\begin{aligned}
1)\,\,\,\,(\,\partial_{a}\Phi_{1}\!\!\mid \,\,=(\,\partial_{\overline{b}}\overline{\Phi_{2} }\! \mid=2x, \,\,\,\,\,\,\,\,\,\,\,\,\,\,2)\,\,\,\,(\,\partial_{a}\Phi_{2}\!\!\mid &=-\,(\,\partial_{\overline{b}}\overline{\Phi_{1} }\! \mid=b 
,\\
3)\,\,\,\,(\,\partial_{a}\Phi_{1}\!\!\mid \,\,=(\,\partial_{b}\Phi_{2}\! \mid=2x, \,\,\,\,\,\,\,\,\,\,\,\,\,\,4)\,\,\,\,(\,\partial_{\overline{a}}\Phi_{2}\!\!\mid &=-\,(\,\partial_{\overline{b}}\Phi_{1}\!\! \mid=b 
\end{aligned}
\right. 
\end{equation*}
in coordinates $x,z,u$ of 3D space. We see that equations  (\ref{Eq3}) are fulfilled, i.e. the 4-dimensional quaternionic function $\psi\left(p\right)=p^{2}$ is H-holomorphic in\! $\mathbb{H}$; equations (\ref{Eq3}-1) and (\ref{Eq3}-3) as well as equations (\ref{Eq3}-2) and (\ref{Eq3}-4) are correspondingly identical, and hence the left quaternionic derivative equals the right one after the transition to 3D space.
\subsection{Example of natural exponential function}  \label{Exa2}
 \textit{The function $\psi\left(p\right)=e^{p}=\Phi_{1}+\Phi_{2}\cdot j,$ where $e$ is the base of the natural logarithm.} We represent the quaternion variable $p=x+yi+zj+uk$ as a sum of real and imaginary parts: $p=x+vr,$ where $v=\sqrt{y^{2}+z^{2}+u^{2}}$ is a real value and 
 \begin{equation}
\label{Eq4} r=\frac{yi+zj+uk}{\sqrt{y^{2}+z^{2}+u^{2}}}
\end{equation}
is a purely imaginary unit quaternion, so its square is $-1.$  We will consider further $v > 0.$ Since $r^{2}=-1$ as well as $x$ and $v$ are real values, the quaternionic formula $p=x+vr$ is algebraically equivalent to the complex formula $\xi=x+yi$ (here and further, if needed, we use the designation $ \xi$ for a complex variable instead of the usual in complex analysis designation $z$ in order to avoid confusing with the quaternionic component $z$). Then, using the quaternionic analogue of Euler's formula: $ e^{vr}=\cos v+r \sin v,$ we have
$$\psi\left(p\right)=\Phi_{1}+\Phi_{2}\cdot j=e^{p}=e^{\left(x+vr\right)}=e^{x}e^{vr}=e^{x}\left(\cos v+r\sin v\right)$$
$$=e^{x}\left(\cos v+\frac{yi\sin v}{v}\right)+e^{x}\frac{\left(z+ui\right)\sin v}{v}\cdot j, $$whence%muster jf intervals
$$\Phi_{1}\!=e^{x}\left(\cos v+\frac{yi\sin v}{v}\right)\,\,\,\,\,\,\,\, \text{and}\,\,\,\,\,\,\,\,\Phi_{2}=e^{x}\frac{\left(z+ui\right)\sin v}{v}.$$

Using the expressions $x=\frac{a+\overline{a}}{2},\, y=\frac{a-\overline{a}}{2i},\,z=\frac{b+\overline{b}}{2},\,u=\frac{b-\overline{b}}{2i},$ following from (\ref{Eq2}), we obtain the expressions for $\Phi_{1}$ and $\Phi_{2}$ as functions of $a,\,\overline{a},\,b,\,\overline{b}:$
$$\Phi_{1}\!=2\beta\cos v+\frac{\beta\left(a-\overline{a}\right)\sin v}{v}\,\,\,\,\,\,\,\, \text{and}\,\,\,\,\,\,\,\,\Phi_{2}=\frac{2\beta\, b\sin v}{v},$$ where
\begin{equation}
\label{Eq5} \begin{gathered}
\beta=\frac{e^{\frac{a+\overline{a}}{2}}}{2}=\overline{\beta}, \,\,\,\,\,\,v=\sqrt{y^{2}+\vert b \vert^{2}}=\frac{\sqrt{\left[4\vert p\vert^{2}-\left(a+\overline{a} \right)^{2}\right]}}{2}=\overline{v},\\ \vert p\vert =\sqrt{x^{2}+y^{2}+z^{2}+u^{2}}=\sqrt{a\overline{a}+b\overline{b}}
\end{gathered}
\end{equation}
are real values. We see that in fact we have obtained the expression for quaternionic function $\psi\left(p\right)=e^{p}$ from the expression for complex function $\psi\left(\xi\right)=e^{\xi}$ by the direct replacement of a complex variable $\xi=x+y i$ by a quaternion variable  $p=x+vr.$ Since the complex imaginary unit $i$ must be replaced by its quaternionic analogue $r$\!, we need to provide for such a replacement everywhere when constructing quaternionic holomorphic functions.

Over 40 various H-holomorphic functions, obtained by the replacement of a complex variable by a quaternionic one in expressions for initial complex holomorphic functions, are represented in \cite{pm:ph}. 

The partial derivatives of $\Phi_{1}$ and $\Phi_{2}$ for the function $e^{p}$ are the following:
\begin{equation*}
\begin{aligned}
&\partial_{a}\Phi_{1}=\beta\left[\cos v+\frac{\left(a-\overline{a}+1\right)\sin v}{v}-\frac{\left(a-\overline{a}\right)^{2}\left(v \cos v-\sin v\right)}{4v^{3}}\right],\\ &\partial_{b}\Phi_{2}=\partial_{\overline{b}}\overline{\Phi_{2}}=\beta\left[\frac{2\sin v}{v}+\frac{b\overline{b}\left(v\cos v-\sin v\right)}{v^{3}}\right],\\ &\partial_{a}\Phi_{2}=-\partial_{\overline{b}}\Phi_{1}=\beta b\left[\frac{\sin v}{v}-\frac{\left(a-\overline{a}\right)\left(v \cos v-\sin v\right)}{2v^{3}}\right],\\ &\partial_{\overline{a}}\Phi_{2}=-\partial_{\overline{b}}\overline{\Phi_{1} }=\beta b\left[\frac{\sin v}{v}+\frac{\left(a-\overline{a}\right)\left(v \cos v-\sin v\right)}{2v^{3}}\right]. 
\end{aligned}
\end{equation*}

Performing the transition $a=\overline{a}=x$ in them and taking into consideration that $v=\vert b\vert >0$ and $\beta=\frac{e^{x}}{2}$ after this transition as well as $b\overline{b}=\vert b \vert^{2}$ we get %muster  3
\begin{equation*}  \left\{
\begin{aligned}
1)\,(\,\partial_{a}\Phi_{1}\!\!\mid \,\,=(\partial_{\overline{b}}\overline{\Phi_{2} }\! \mid=\frac{e^{x}\left(\cos\vert b \vert+\vert b \vert^{-1}\sin \vert b \vert\right)}{2}\!; \,\,2)\,(\,\partial_{a}\Phi_{2}\!\!\mid &=-\,(\,\partial_{\overline{b}}\overline{\Phi_{1} }\! \mid=\frac{e^{x}b\vert b \vert^{-1}\sin \vert b \vert}{2} 
\!;\\
3)\,(\,\partial_{a}\Phi_{1}\!\!\mid \,\,=(\partial_{b}\Phi_{2} \! \mid=\frac{e^{x}\left(\cos\vert b \vert+\vert b \vert^{-1}\sin \vert b \vert\right)}{2}\!;\,\,4)\,(\,\partial_{\overline{a}}\Phi_{2}\!\!\mid &=-\,(\,\partial_{\overline{b}}\Phi_{1}\!\! \mid=\frac{e^{x}b\vert b \vert^{-1}\sin \vert b \vert}{2}\! .
\end{aligned}
\right. 
\end{equation*}
We see that equations  (\ref{Eq3}) are fulfilled and the function $\psi\left(p\right)=e^{p}$ is H-holomorphic everywhere in $\mathbb{H}$, including the point $p=0$. Indeed, at this point we have equations (\ref{Eq3}) as follows: 
\begin{equation*}  \left\{
\begin{aligned}
1)\,(\,\partial_{a}\Phi_{1}\!\!\mid \,\,=(\partial_{\overline{b}}\overline{\Phi_{2} }\! \mid=1; \,\,2)\,(\,\partial_{a}\Phi_{2}\!\!\mid &=-\,(\,\partial_{\overline{b}}\overline{\Phi_{1} }\!\! \mid=0
;\\
3)\,(\,\partial_{a}\Phi_{1}\!\!\mid \,\,=(\partial_{b}\Phi_{2} \! \mid=1;\,\,4)\,(\,\partial_{\overline{a}}\Phi_{2}\!\!\mid &=-\,(\,\partial_{\overline{b}}\Phi_{1}\!\! \mid=0 ,
\end{aligned}
\right. 
\end{equation*}
which are also fulfilled.

Equations 1) and 3) as well as 2) and 4) are respectively identical and hence the left quaternionic derivative is equal to the right one after the transition $a=\overline{a}=x.$
\section{The full quaternionic derivative}
It was established in \cite{pm:aqg} that the quaternionic generalization of the complex derivative has the following expression for the full (\textit{uniting the parts of left and right derivatives} \cite[\!p.~18]{pm:oph}) quaternionic derivative  of the $k'\text{th}$ order:
\begin{equation}
\label{Eq6} \psi_H^{\left(k\right)}\!\!\left(p\right)=\Phi_1^{(k)}+\Phi_2^{(k)}\cdot j,  %asdfghjklpoiuztrewqa  asdfghj 
\end{equation}
where the constituents $\Phi_1^{(k)}$ and $\Phi_2^{(k)}$ are expressed by
$$\Phi_1^{(k)}=\partial_{a}\Phi_1^{(k-1)}\!\!+\partial_{\overline{a}}\Phi_1^{(k-1)} \,\,\,\text{and}\,\,\,\,\,\,\,\Phi_2^{(k)}=\partial_{a}\Phi_2^{(k-1)}\!\!+\partial_{\overline{a}}\Phi_2^{(k-1)},$$ and
$\Phi_1^{(k-1)} \text{and}\,\, \Phi_2^{(k-1)}$ are the constituents of the $\left(k-1\right)^{'}\!\text{th}$ full derivative of  $\psi_{H}\left(p\right)$ represented in the Cayley–Dickson doubling form as $ \psi_H^{\left(k-1\right)}\!\!\left(p\right)=\Phi_1^{(k-1)}+\Phi_2^{(k-1)}\cdot j,  k\geq1;$ $\Phi_1^{(0)}=\Phi_{1}\!\!\left(a,b,\overline{a}, \overline{b} \right) \text{and}\,\, \Phi_2^{(0)}=\Phi_{2}\!\!\left(a,b,\overline{a}, \overline{b} \right).$

 \textit{If a quaternion function $\psi\left(p\right)$ is once H-differentiable in $G_{4}\subset \mathbb{H},$ then it possesses the full quaternionic derivatives of all orders in $G_{4}$ (except, possibly, at certain singularities), each one H-differentiable} \cite{pm:aqg}. 
Analogously to complex analysis there are various quaternionic equivalents  of formula (\ref{Eq6}). See details in \cite{pm:aqg,pm:oph}.
 
 The full first quaternionic derivative of the function $\psi_{H}\left(p\right)=p^{2}$ (see Example  \ref{Exa1} ) is the following:
 \begin{equation*}
\begin{aligned}
 \psi_H^{(1)}\!\!\left(p\right)&=\left(p^{2}\right)^{{\left(1\right)}}=\Phi_1^{(1)}+\Phi_2^{(1)}\cdot j=\left(\partial_{a}\Phi_1\!+\partial_{\overline{a}}\Phi_1\right)+\left(\partial_{a}\Phi_2\!+\partial_{\overline{a}}\Phi_2\right)\cdot j \\
& =\left(2a+0\right)+\left(b+b\right)\cdot j=2a+2b \cdot j=2p.
\end{aligned}
\end{equation*}
Given that the function $e^{p}$ has  $\partial_{\overline{a}}\Phi_{1}=\beta\left[\frac{v \cos v-\sin v}{v}+\frac{\left(a-\overline{a}\right)^{2}\left(v \cos v-\sin v\right)}{4v^{3}}\right],$ we obtain after some computation the expression for the first full quaternionic derivative of the function $\psi_{H}\left(p\right)=e^{p}$ (see Example \ref{Exa2}) as follows:
\begin{equation*}
 \begin{split}
 \psi_H^{(1)}\!\!\left(p\right)&=\left(e^{p}\right)^{{\left(1\right)}}=\Phi_1^{(1)}+\Phi_2^{(1)}\cdot j=\left(\partial_{a}\Phi_1\!+\partial_{\overline{a}}\Phi_1\right)+\left(\partial_{a}\Phi_2\!+\partial_{\overline{a}}\Phi_2\right)
 \cdot j \\
 & = \left\{ \beta\left[\cos v+\frac{\left(a-\overline{a}+1\right)\sin v}{v}-\frac{\left(a-\overline{a}\right)^{2}\left(v \cos v-\sin v\right)}{4v^{3}}\right]  \right.  
\\ &  \left.  \,\,\,\,\,\,\,\,\,\,\,\, +\beta\left[\frac{v \cos v-\sin v}{v} +\frac{\left(a-\overline{a}\right)^{2}\left(v \cos v-\sin v\right)}{4v^{3}}\right] \right\} \\  & 
 +\left\{\beta b\left[\frac{\sin v}{v}-\frac{\left(a-\overline{a}\right)\left(v \cos v-\sin v\right)}{2v^{3}}\right]  \right. \\
  &\left.  \,\,\,\,\,\,\,\,\,\,\,\,  +\beta b\left[\frac{\sin v}{v}+\frac{\left(a-\overline{a}\right)\left(v \cos v-\sin v\right)}{2v^{3}}\right]\right\}\cdot j\\
 & =2\beta\cos v+\frac{\beta\left(a-\overline{a}\right)\sin v}{v}+\frac{2\,\beta\, b\sin v}{v}\cdot j=e^{p}.
    \end{split}
\end{equation*}

We see that the expressions for the full first quaternionic derivatives of the H-~holomorphic functions  $\psi_{H}\!\left(p\right)=p^{2}$ and $\psi_{H}\!\left(p\right)=e^{p}$ have the same forms as in real and complex analysis. The same there is for the derivatives of higher orders. We also see that the full quaternionic derivatives of all orders of the natural exponential function are equal to the function itself just as in complex (and real) analysis.
\section{The new class of H-holomorphic functions}

We introduce a new class of H-holomorphic functions, which consists of the quaternionic functions meeting equations (\ref{Eq3}). In accordance with \cite{pm:eac} we point out here some properties of this class :

1) Each H-holomorphic function $\psi_{H}\left(p\right)$ can be constructed from its complex holomorphic analogue $\psi_{C}\left(\xi\right)$ (without change of a functional dependence form) by replacing a complex variable $\xi\in G_{2}\subseteq \mathbb{C}$ as a single whole in an expression for $\psi_{C}\left(\xi\right)$ by a quaternionic variable $p\in G_{4}\subseteq \mathbb{H},$ where $G_{4}$ is defined (except, possibly, at certain singularities) by the relation $G_{4}\supset G_{2}$ in the sense that $G_{2}$ exactly follows from $G_{4}$ upon the transition from $p$ to $\xi$ \cite[p.~15]{pm:eac}.  At that we must first convert the initial expression for a complex holomorphic function into the expression in which there is dependence on complex variable (and/or its conjugation) only as a single whole.

2) The constituents of the H-holomorphic functions (and their derivatives of all orders) in the Cayley-Dickson doubling form (\ref{Eq1}) have the following general representation forms \cite{pm:oph}: $$\Phi_{1}\!\!\left(a,b,\overline{a}, \overline{b} \right)=A\left[a,\overline{a},\left(b \overline{b}\right)\!\right]\!, \,\Phi_{2}\!\!\left(a,b,\overline{a},\overline{b} \right)=B\!\left[\left(a \overline{a}\right)\!,\left(a \overline{a}\right)_{m}\!,\left(b \overline{b}\right)\right]\!\!b=Bb, \,B=\overline{B},$$ where $\left(a \overline{a}\right)_{m}$ denotes the form $a^{m}\overline{a}^{0}+a^{{\left(m-1\right)}}\overline{a}^{1}+a^{{\left(m-2\right)}}\overline{a}^{2}+\cdots+a^{0}\overline{a}^{m}, m=0,1,2,3, \dots $ with the number of summands $\left(m+1\right).$ Together with this form or instead of it, there can be another symmetric forms that are invariant under complex conjugation \cite{pm:oph,pm:aqg}. Such the forms of $\Phi_{1}$ and $\Phi_{2}$ are typical of H-holomorphic functions and their derivatives. They can serve as gauge for correctness of results obtained when constructing H-holomorphic functions. Such a symmetry in cases of derivatives is a consequence of uniting unsymmetrical parts of the left and right derivatives \cite{pm:oph}. All examples in this article illustrate such the forms (see also \cite{pm:ph}) for H-holomorphic functions and their full quaternionic derivatives.

3) Algebraic properties of this class are all identical to ones of the class of complex holomorphic functions.\! In particular, the quaternionic multiplication of the H-holomorphic functions behaves as commutative and the left quotient of two H-holomorphic functions is equal to the right one \cite{pm:eac}.

4) The formulae for differentiating sums, products, ratios, inverses, and compositions of H-holomorphic functions are all identical to their complex and real analogues \cite{pm:eac}. The expressions for the full quaternionic derivatives of the H-holomorphic functions are all identical to the expressions for corresponding derivatives of complex holomorphic analogues. For example, if the first derivative of the complex holomorphic function $\psi_{C}\left(\xi\right)=\xi^{3}$ is $3\xi^{2},$ then the first full quaternionic derivative of the H-holomorphic analogue $\psi_{H}\left(p\right)=p^{3}$ is $3p^{2}$ \cite{pm:aqg,pm:eac}. One can just verify these properties, constructing H-holomorphic functions from their complex holomorphic analogues. The H-holomorphic functions possess the full H-holomorphic quaternionic derivatives of all orders. We can denote the presented class of  H-holomorphic functions by $\text{H}^{\infty}.$  Note that the known Cauchy-Riemann equations  \cite{mh:ca} follow from their quaternionic generalization (\ref{Eq3}) upon the transition to complex plane  \cite[\!p.~29]{pm:aqg}.

5) The presented concept of quaternionic holomorphy allows us to explore steady state vector fields in 3D space, each of them corresponds to some H-holomorphic function considered as a quaternionic potential \cite{pm:qpc}.
\subsection{Example of reciprocal function}  \label{Exa3}
 \textit{The reciprocal function} $\psi\left(p\right)=p^{{-1}},\, \,\,p\neq 0.$ The corresponding complex holomorphic analogue is $\psi\left(\xi\right)=\xi^{{-1}}=\frac{\overline{\xi}}{\vert \xi\vert^{2} },\, \,\,\xi\neq 0.$ By replacing a complex variable $\xi$ by a quaternionic $p$ in it we obtain the following expression for the H-holomorphic function: $\psi\left(p\right)=p^{{-1}}=\frac{\overline{p}}{\vert p\vert^{2}}=\frac{\overline{a}}{\left(a \overline{a}+b\overline{b}\right)}-\frac{b}{\left(a \overline{a}+b\overline{b}\right)}\cdot j=\Phi_{1}+\Phi_{2}\cdot j,\,\,\,\, p\neq0,$ whence $\Phi_{1}=\frac{\overline{a}}{\left(a \overline{a}+b\overline{b}\right) }$ and $\Phi_{2}=-\frac{b}{\left(a \overline{a}+b\overline{b}\right) }.$ Correspondingly, $\overline{\Phi_{1}}=\frac{a}{\left(a \overline{a}+b\overline{b}\right) }$ and $\overline{\Phi_{2}}=-\frac{\overline{b}}{\left(a \overline{a}+b\overline{b}\right)}.$ 

To check the H-holomorphy of the function $\psi\left(p\right)=p^{{-1}}$ we compute the partial derivatives of the functions $\Phi_{1}$ and $\Phi_{2}\!\!: \partial_{a}\Phi_{1}=-\frac{\overline{a}^{2}}{\left(a \overline{a}+b\overline{b}\right)^{2}}, \,\,\, \partial_{\overline{b}}\overline{\Phi_{2} }=\partial_{b}\Phi_{2}=-\frac{a\overline{a}}{\left(a \overline{a}+b\overline{b}\right)^{2}},\,$ and $\partial_{a}\Phi_{2}=-\partial_{\overline{b}}\Phi_{1}=\frac{\overline{a}b}{\left(a \overline{a}+b\overline{b}\right)^{2}},\,\,\, \partial_{\overline a}\Phi_{2}=-\partial_{\overline{b}}\overline{\Phi_{1}}=\frac{ab}{\left(a \overline{a}+b\overline{b}\right)^{2}}.$ After performing the transition $ a=\overline{a}=x$ we see that H-holomorphy equations (\ref{Eq3}) are fulfilled:
\begin{equation*}  \left\{
\begin{aligned}
1)\,\,\,(\,\partial_{a}\Phi_{1}\!\!\mid \,\,=(\,\partial_{\overline{b}}\overline{\Phi_{2} }\! \mid=-\frac{x^{2}}{\left(x^{2}+b\overline{b}\right)^{2}}, \,\,\,\,\,\,\,\,\,2)\,\,\,(\,\partial_{a}\Phi_{2}\!\!\mid &=-\,(\,\partial_{\overline{b}}\overline{\Phi_{1} }\! \mid=\frac{xb}{\left(x^{2}+b\overline{b}\right)^{2}} 
,\\
3)\,\,\,(\,\partial_{a}\Phi_{1}\!\!\mid \,\,=(\,\partial_{b}\Phi_{2}\! \mid=-\frac{x^{2}}{\left(x^{2}+b\overline{b}\right)^{2}}, \,\,\,\,\,\,\,\,\,4)\,\,\,(\,\partial_{\overline{a}}\Phi_{2}\!\!\mid &=-\,(\,\partial_{\overline{b}}\Phi_{1}\!\! \mid=\frac{xb}{\left(x^{2}+b\overline{b}\right)^{2}} .
\end{aligned}
\right. 
\end{equation*}\\The function $\psi\left(p\right)=p^{{-1}}$ is H-holomorphic at the points $p\in\!  \mathbb{H}\! \setminus \!\! \left\{0\right\}$ and its left and right quaternionic derivatives become equal after the transition $ a=\overline{a}=x$ to 3D space. \\
\textit{The first derivative of the function $\psi_{H}\!\left(p\right)=p^{{-1}}$\!.} Using the expressions for $\partial_{a}\Phi_{1},\,\partial_{a}\Phi_{2}, \\ \partial_{\overline{a}}\Phi_{2}$ as well as $\partial_{\overline{a}}\Phi_{1}=\frac{b\overline{b}}{\left(a\overline{a}+b\overline{b}\right)^{2}}$ in (\ref{Eq6}), we get the following expression for the first full (further we omit, for brevity,  the word “\,full”) quaternionic derivative of the function $\psi_{H}\!\left(p\right)=p^{{-1}}\!:\\$
\begin{align*}
&\left(p^{{-1}}\right)^{{\left(1\right)}}=\Phi_1^{(1)}+\Phi_2^{(1)}\cdot j=\left(\partial_{a}\Phi_1\!+\partial_{\overline{a}}\Phi_1\right)+\left(\partial_{a}\Phi_2\!+\partial_{\overline{a}}\Phi_2\right)
 \cdot j \\ &=\left[-\frac{\overline{a}^{2}}{\left(a \overline{a}+b\overline{b}\right) ^{2}}+\frac{b\overline{b}}{\left(a \overline{a}+b\overline{b}\right) ^{2}}\right]+\left[\frac{\overline{a}\,b}{\left(a \overline{a}+b\overline{b}\right) ^{2}}+\frac{a\,b}{\left(a \overline{a}+b\overline{b}\right) ^{2}}\right]\cdot j\\&=-\frac{\left[\overline{a}^{2}-b \overline{b}-b\left(a+\overline{a}\right)\cdot j\right]}{\left(a \overline{a}+b\overline{b}\right) ^{2}}=-\frac{\overline{p}^{2}}{\left(a \overline{a}+b\overline{b}\right) ^{2}}=-\frac{\overline{p}^{2}}{\vert p\vert^{4}}.
\end{align*}
Given $p^{{-2}}\!=p^{{-1}}\cdot p^{{-1}}=\frac{\overline{p}^{2}}{\vert p\vert^{4}},$ we get the following expression for the first derivative of the function $\psi_{H}\left(p\right)=p^{{-1}}\!:$
\begin{equation} \label{Eq7}
\left(p^{{-1}}\right)^{{\left(1\right)}}=\Phi_1^{(1)}+\Phi_2^{(1)}\cdot j=-p^{{-2}},
\end{equation}
where $\Phi_1^{(1)}=\frac{b\overline{b}-\overline{a}^{2}}{{\left(a \overline{a}+b\overline{b}\right) ^{2}}},\,\Phi_2^{(1)}=\frac{b\left(a+\overline{a}\right)}{{\left(a \overline{a}+b\overline{b}\right) ^{2}}}.$ Their complex conjugations are $\overline{\Phi_1^{(1)}}=\frac{\overline{b}b -a^{2}}{{\left(a \overline{a}+b\overline{b}\right) ^{2}}}, \\\overline{\Phi_2^{(1)}}=\frac{\overline{b}\left(\overline{a}+a\right)}{{\left(a \overline{a}+b\overline{b}\right) ^{2}}}.$ We see that the obtained expression (\ref{Eq7}) for the first quaternionic derivative of the function $p^{{-1}}$ is similar to the expression for the complex analogue: $\left(\xi^{{-1}}\right)^{{\left(1\right)}}=-\xi^{{-2}}.$

To verify the H-holomorphy of the first derivative $\left(p^{{-1}}\right)^{{\left(1\right)}}$ we calculate the following partial derivatives: $\partial_{a}\Phi_1^{(1)}=\frac{2\overline{a}\left(\overline{a}^{2}-b\overline{b}\right)}{\left(a \overline{a}+b\overline{b}\right) ^{3}},\,\partial_{b}\Phi_2^{(1)}=\partial_{\overline{b}}\overline{\Phi_2^{(1)}}=\frac{\left(a+\overline{a}\right)\left(a\overline{a}-b\overline{b}\right)}{\left(a \overline{a}+b\overline{b}\right) ^{3}}$ and $\partial_{a}\Phi_2^{(1)}=-\partial_{\overline{b}}\Phi_1^{(1)}= \frac{b\left(b\overline{b}-a\overline{a}-2\overline{a}^{2}\right)}{\left(a \overline{a}+b\overline{b}\right) ^{3}},\,\partial_{\overline{a}}\Phi_2^{(1)}=-\partial_{\overline{b}}\overline{\Phi_1^{(1)}}=   \frac{b\left(b\overline{b}-a\overline{a}-2a^{2}\right)}{\left(a \overline{a}+b\overline{b}\right) ^{3}}.$ After performing the transition  $ a=\overline{a}=x$ we see that the H-holomorphy equations (\ref{Eq3}) are fulfilled for the first derivative of the function $\psi_{H}\!\left(p\right)=p^{{-1}}:$
%muster for long eq (3)
\begin{equation*}  \left\{ 
\begin{aligned}
1)\,\,\,\,(\,\partial_{a}\Phi_1^{(1)}\!\!\mid \,\,=(\,\partial_{\overline{b}}\overline{\Phi_{2} ^{(1)}}\!\!  \mid=  \frac{2x\left(x^{2}-b\overline{b}\right)}{\left(x^{2}+b\overline{b}\right)^{3}}, \,\,\,\,\,\,\,\,\,\,2)\,\,\,\,(\,\partial_{a}\Phi_{2}^{(1)}\!\!\mid &=-\,(\,\partial_{\overline{b}}\overline{\Phi_{1}^{(1)}}\!\! \mid    =\frac{b\left(b\overline{b}-3x^{2}\right)}{\left(x^{2}+b\overline{b}\right)^{3}},\\
3)\,\,\,\,(\,\partial_{a}\Phi_1^{(1)}\!\!\mid \,\,=(\,\partial_{b}\Phi_2^{(1)}\!\! \mid=\frac{2x\left(x^{2}-b\overline{b}\right)}{\left(x^{2}+b\overline{b}\right)^{3}}, \,\,\,\,\,\,\,\,\,\,4)\,\,\,\,(\,\partial_{\overline{a}}\Phi_2^{(1)}\!\!\mid &=-\,(\,\partial_{\overline{b}}\Phi_1^{(1)}\!\! \mid=     \frac{b\left(b\overline{b}-3x^{2}\right)}{\left(x^{2}+b\overline{b}\right)^{3}}  .
\end{aligned}
\right. 
\end{equation*}

The first derivative of the function $\psi_{H}\!\left(p\right)=p^{{-1}}$ is also H-holomorphic at  $p\in \mathbb{H}\setminus \!\left\{0\right\}$ and its left and right quaternionic derivatives become equal after the transition $a=\overline{a}=x$ to 3D space. \\
\textit{The second derivative of the function} $\psi_{H}\!\left(p\right)=p^{{-1}}.$ Quite analogously, using formula (\ref{Eq6}) for $k=2$ and the calculated derivatives $\partial_{a}\Phi_1^{(1)}\!,\,\,\partial_{a}\Phi_2^{(1)}\!,\,\, \partial_{\overline{a}}\Phi_2^{(1)}$ as well as $\partial_{\overline{a}}\Phi_1^{(1)}=-\frac{2b\overline{b}\left(a+\overline{a}\right)}{\left(a \overline{a}+b\overline{b}\right) ^{3}},$ we obtain after some calculations the following expression for the second derivative of the function $\psi_{H}\!\left(p\right)=p^{{-1}}:$
 \begin{align*} 
 \left(p^{{-1}}\right)^{{\left(2\right)}} &=\Phi_1^{(2)}+\Phi_2^{(2)}\cdot j=   \left(\partial_{a}\Phi_1^{(1)}\!+\partial_{\overline{a}}\Phi_1^{(1)}\right)+\left(\partial_{a}\Phi_2^{(1)}\!+\partial_{\overline{a}}\Phi_2^{(1)}\right)
 \cdot j \\ &= \frac{2\left[\overline{a}^{3}-2\overline{a}b\overline{b}-ab\overline{b}\right]+2\left[b\overline{b}-\overline{a}^{2}-a\left(\overline{a}+a\right)\right]b\cdot j}{\left(a\overline{a}+b\overline{b}\right)^{3}},
\end{align*}
 whence $$\Phi_1^{(2)}=\frac{2\left[\overline{a}^{3}-\left(2\overline{a}+a\right)b\overline{b}\right]}{\left(a \overline{a}+b\overline{b}\right) ^{3}},\,\,\,\Phi_2^{(2)}=\frac{2\left[b\overline{b}-\left(a^{2}+a\overline{a}+\overline{a}^{2}\right)\right]b}{\left(a \overline{a}+b\overline{b}\right) ^{3}}=\frac{2\left[b\overline{b}-\left(a\overline{a}\right)_{2}\right]b}{\left(a \overline{a}+b\overline{b}\right) ^{3}}.$$
 
Taking into account 
\begin{align*}
\overline{p}^{2} &=\left(\overline{a}-b\cdot j\right)\cdot \left(\overline{a}-b\cdot j\right)=\overline{a}^{2}-b\overline{b}-b\left(a+\overline{a}\right)\!\cdot j,\\      \overline{p}^{3}&=\overline{p}^{2}\cdot\overline{p}=\left[\overline{a}^{2}-b\overline{b}-b\left(a+\overline{a}\right)\!\cdot j\right]\cdot\left(\overline{a}-b\cdot j\right) \\ =&\left[\overline{a}^{3}-2\overline{a}\, b\overline{b}-ab\overline{b}\right]+\left[b\overline{b}-\overline{a}^{2}-a\left(\overline{a}+a\right)\right]\!b\cdot j,  \\ \,\, \,\,\,\,\,\,\,\,\,\,\,\,&p^{{-3}}=p^{{-1}}\!\cdot p^{{-1}}\!\cdot  p^{{-1}}=\frac{\overline{p}^{3}}{\vert p\vert^{6}}=\frac{\overline{p}^{3}}{\left(a \overline{a}+b\overline{b}\right)^{3}},
\end{align*}
we get the following expression for the second derivative of the function $\psi_{H}\!\left(p\right)=p^{{-1}}:$
$$\left(p^{{-1}}\right)^{{\left(2\right)}}=2\frac{\overline{p}^{3}}{\left(a \overline{a}+b\overline{b}\right)^{3}}=2p^{{-3}}.$$
We see that the obtained expression has the same form as the complex (or real) analogue: $\left(\xi^{{-1}}\right)^{{\left(2\right)}}=2\xi^{{-3}}.$
To verify the H-holomorphy of the second derivative $\left(p^{{-1}}\right)^{{\left(2\right)}}$ we calculate the following derivatives:
\begin{align*} %Muster align ohne Number
       &\partial_{a}\Phi_1^{(2)}=\frac{2\left[-3\overline{a}^{4}+2\overline{a}\left(3\overline{a}+a\right)b\overline{b}\,-b^{2}\overline{b}^{2}\right]}{\left(a \overline{a}+b\overline{b}\right)^{4}},     \\ &\partial_{b}\Phi_2^{(2)}=\partial_{\overline{b}}\overline{\Phi_2^{(2)}} =\frac{2\left[-a\overline{a}\left(a\overline{a}\right)_{2}+2\left(a^{2}+2a\overline{a}+\overline{a}^{2}\right)b\overline{b}-b^{2}\overline{b}^{2}\right]}{\left(a \overline{a}+b\overline{b}\right) ^{4}},  \\ &\partial_{a}\Phi_2^{(2)} =-\partial_{\overline{b}}\Phi_1^{(2)}=\frac{  2\overline{a}  \left(a^{2}+2a\overline{a}+3\overline{a}^{2}\right)b}{\left(a \overline{a}+b\overline{b}\right) ^{4}}-\frac{4\left(a+2\overline{a}\right)b^{2}\overline{b}}{\left(a \overline{a}+b\overline{b}\right) ^{4}}, \\ &   \partial_{\overline {a}}\Phi_2^{(2)}=-\partial_{\overline{b}}\overline{\Phi_1^{(2)}}= \frac{2\left[a\overline{a}\left(\overline{a}+2a\right)b+3a^{3}b-2\left(\overline{a}+2a\right)b^{2}\overline{b}\right]}{\left(a \overline{a}+b\overline{b}\right) ^{4}}.            
\end{align*}
After performing the transition $ a=\overline{a}=x$ to 3D space we see that the H-holomorphy equations (\ref{Eq3}) are fulfilled for the second derivative of the function $\psi_{H}\!\left(p\right)=p^{{-1}}\!:$
\begin{equation*}  \left\{ 
\begin{aligned}
&1)\,\,\,\,(\,\partial_{a}\Phi_1^{(2)}\!\!\mid \,\,=(\,\partial_{\overline{b}}\overline{\Phi_{2} ^{(2)}}\!\!  \mid =  \frac{2\left(-3x^{4}+8x^{2}b\overline{b}-b^{2}\overline{b}^{2}\right)}{\left(x^{2}+b\overline{b}\right) ^{4}}, \\ &2)\,\,\,\,(\,\partial_{a}\Phi_{2}^{(2)}\!\!\mid\, =-\,(\,\partial_{\overline{b}}\overline{\Phi_{1}^{(2)}}\!\! \mid    =\frac{12xb\left(x^{2}-b\overline{b}\right)}{\left(x^{2}+b\overline{b}\right) ^{4}},\\
&3)\,\,\,\,(\,\partial_{a}\Phi_1^{(2)}\!\!\mid \,\,=(\,\partial_{b}\Phi_2^{(2)}\!\! \mid=\frac{2\left(-3x^{4}+8x^{2}b\overline{b}-b^{2}\overline{b}^{2}\right)}{\left(x^{2}+b\overline{b}\right) ^{4}}, \,\,\,\,\,\,\,\, \\ &4)\,\,\,\,(\,\partial_{\overline{a}}\Phi_2^{(2)}\!\!\mid \,=-\,(\,\partial_{\overline{b}}\Phi_1^{(2)}\!\! \mid=     \frac{12xb\left(x^{2}-b\overline{b}\right)}{\left(x^{2}+b\overline{b}\right) ^{4}}.
\end{aligned}
\right. 
\end{equation*}

The second derivative of the H-holomorphic function $\psi_{H}\!\left(p\right)=p^{{-1}}$ is H-holomorphic at the points $p\in \mathbb{H}\setminus \!\left\{0\right\}$ too and its left and right quaternionic derivatives become equal after the transition $a=\overline{a}=x$ to 3D space. 
%Similarly one can verify that the derivatives of higher orders are H-holomorphic.

Already this example shows that calculations with quaternionic functions can be cumbersome, error-prone and tedious manual procedures. To avoid this problem there is a special programmes pack \cite{pm:ph} in the Wolfram Mathematica programming language, which verifies the commutative behavior of multiplication of H-holomorphic functions, calculates their quaternionic derivatives of all orders, verifies whether or not the quaternionic functions (or their derivatives) satisfy quaternionic generalization (\ref{Eq3}) of Cauchy-Riemann's equations, calculates quaternionic expressions of any complexity as well as expressions for 3D potential fields, etc.  Using this programmes pack we can  easily verify that the derivatives of $\psi_{H}\!\left(p\right)=p^{{-1}}$ of higher orders are also all H-holomorphic.
\subsection{Example of natural logarithmic function}  \label{Exa4} 
 \textit{The principal value of the quaternionic natural logarithmic function:} $\psi\!\left(p\right)= \text{Ln}\!\left(p\right).$ The principal value of the initial complex logarithmic function is given \cite{mh:ca} by $\text{Ln}\left(\xi\right)=\ln\vert\xi\vert+i\, \text{arg} \left(\xi\right)=\ln\vert\xi\vert+i\arccos\frac{\xi+\overline{\xi}}{2\vert\xi\vert}$, where $-\pi <\text{arg}\left(\xi\right)\leq\pi$ and $\vert \xi\vert >0$. By replacing a complex variable $\xi$ by by a quaternionic $p$ and $i$ by $r$ (as noted above) we obtain the following expression for the principal value of the quaternionic natural logarithmic function: $$\psi\!\left(p\right) =\text{Ln}\left(p\right)=\ln\vert p\vert+r\cdot\text{Arccos}\frac{p+\overline{p}}{2\vert p\vert},$$ 
 where $\vert p\vert >0$ is defined by  (\ref{Eq5}), $r$ by (\ref{Eq4}) and $-\pi <\text{Arccos}\frac{p+\overline{p}}{2 \vert p\vert}\leq\pi .$ Using (\ref{Eq2}), we finally obtain the following expression for the quaternionic logarithmic function as a function of the variables $a,\,\overline{a},\,b,\,\overline{b}:$
 $$\psi\!\left(p\right) =\text{Ln}\left(p\right)=\Phi_{1}+\Phi_{2}\cdot j=\ln\vert p\vert+\frac{\left(a-\overline{a}\right)\text{Arccos}\frac{a+\overline{a}}{2\vert p\vert}}{2v}+\frac{b \,\text{Arccos}\frac{a+\overline{a}}{2\vert p\vert}}{v}\cdot j,$$\\where $v$ is the same as in (\ref{Eq5}), whence $$\Phi_{1}=\ln\vert p\vert+\frac{\left(a-\overline{a}\right)\text{Arccos}\frac{a+\overline{a}}{2\vert p\vert}}{2v},\,\,\,\Phi_{2}=\frac{b \,\text{Arccos}\frac{a+\overline{a}}{2\vert p\vert}}{v},\,\,\,p,v\neq0.$$
 
 After the computation of the partial derivatives we have
\begin{align*} %Muster align ohne Number
&\partial_{a}\Phi_{1} =\frac{\overline{a}}{2\vert p\vert^{2}}+\frac{\theta}{2v}-\frac{\left(a-\overline{a}\right)\left[2\vert p\vert^{2}-\left(a+\overline{a}\right)\overline{a} \right]}{2\vert p\vert^{2}\left[4 \vert p\vert^{2}-\left(a+\overline{a}\right)^{2}\right]}+\frac{\left(a-\overline{a}\right)^{2}\theta}{\left(2v\right)^{3}},  \\  & \partial_{b}\Phi_{2}=\partial_{\overline{b}}\overline{\Phi_{2}}=\theta\left(\frac{1}{v}-\frac{\vert b\vert^{2} }{2v^{3}}\right)+\frac{\left(a+\overline{a}\right)\vert b\vert^{2} }{\vert p\vert^{2}\left[4 \vert p\vert^{2}-\left(a+\overline{a}\right)^{2}\right]},\\ & \partial_{a}\Phi_{2}=-\partial_{\overline{b}}\Phi_{1}=-\frac{b \left[2\vert p\vert^{2}-\left(a+\overline{a}\right)\overline{a} \right]}{\vert p\vert^{2}\left[4 \vert p\vert^{2}-\left(a+\overline{a}\right)^{2}\right]}+\frac{\left(a-\overline{a}\right)b\theta}{4v^{3}},  \\  &\partial_{\overline {a}}\Phi_{2}=-\partial_{\overline{b}}\overline{\Phi_{1}}=-\frac{b \left[2\vert p\vert^{2}-\left(a+\overline{a}\right)a \right]}{\vert p\vert^{2}\left[4 \vert p\vert^{2}-\left(a+\overline{a}\right)^{2}\right]}-\frac{\left(a-\overline{a}\right)b\theta}{4v^{3}},  
\end{align*}
where $ \theta= \text{Arccos}\frac{a+\overline{a}}{2\vert p\vert}$ and $a,b\neq0.$ After performing the transition $a=\overline{a}=x,$ we use the relations $\vert p\vert=\vert p_{3}\vert=\sqrt{x^{2}+z^{2}+u^{2}}=\sqrt{x^{2}+\vert b\vert^{2}}>0,\,\,\, 4\vert p\vert^{2}-\left(a+\overline{a}\right)^{2}=4\vert p_{3}\vert^{2}-4x^{2}=4\vert b\vert^{2} >0, \,\,\,v=\sqrt{b\overline{b}}=\vert b\vert >0, \,\,\,\theta= \text{Arccos}\frac{x}{\vert p_{3}\vert}.$ Then, substituting these relations into the above expressions for the partial derivatives, we see that H-holomorphy equations (\ref{Eq3}) are fulfilled:
\begin{equation*}  \left\{ 
\begin{aligned}
1)\,\,\,(\,\partial_{a}\Phi_1\!\!\mid \,\,=(\,\partial_{\overline{b}}\overline{\Phi_{2}}\!\!  \mid=  \frac{x}{2 \vert p_{3}\vert^{2}}+\frac{ \text{Arccos}\frac{x}{\vert p_{3}\vert}}{2\vert b\vert}, \,\,\,\,\,\,\,\,\,\,2)\,\,\,(\,\partial_{a}\Phi_{2}\!\!\mid &=-\,(\,\partial_{\overline{b}}\overline{\Phi_{1}}\!\! \mid  =  -\frac{b}{2 \vert p_{3}\vert^{2}},\\
3)\,\,\,(\,\partial_{a}\Phi_1\!\!\mid \,\,=(\,\partial_{b}\Phi_2\!\! \mid=\frac{x}{2 \vert p_{3}\vert^{2}}+\frac{ \text{Arccos}\frac{x}{\vert p_{3}\vert}}{2\vert b\vert}, \,\,\,\,\,\,\,\,\,\,4)\,\,\,(\,\partial_{\overline{a}}\Phi_2\!\!\mid &=-\,(\,\partial_{\overline{b}}\Phi_1\!\! \mid= -\frac{b}{2 \vert p_{3}\vert^{2}}.
\end{aligned}
\right. 
\end{equation*}

The obtained quaternionic analogue $\psi\!\left(p\right)= \text{Ln}\!\left(p\right)$ of the principal value of complex natural logarithmic function is H-holomorphic in $\mathbb{H}\setminus \!\left\{0\right\}$. The left and the right quaternionic derivatives become equal after the transition to 3D space.\\ \textit{The first derivative of the function $\psi_{H}\left(p\right)=\text{Ln}\!\left(p\right)$}. Taking into consideration the fact that $\partial_{\overline{a}}\Phi_1 = \frac{a}{2\vert p\vert^{2}}-\frac{\theta}{2v}-\frac{\left(a-\overline{a}\right)^{2}\theta}{\left(2v\right)^{3}}-\frac{\left(a-\overline{a}\right)\left[2\vert p\vert^{2}-\left(a+\overline{a}\right)a \right]}{2\vert p\vert^{2}\left[4 \vert p\vert^{2}-\left(a+\overline{a}\right)^{2}\right] },$ we get the following expression for the first derivative of the function 
$\psi_{H}\!\left(p\right)=\text{Ln}\left(p\right):$
\begin{align*} %Muster c & ohne Number
 \left(\text{Ln}\left(p\right)\right)^{{\left(1\right)}}&=\Phi_1^{(1)}+\Phi_2^{(1)}\cdot j=
\left(\partial_{a}\Phi_1\!+\partial_{\overline{a}}\Phi_1\right)+\left(\partial_{a}\Phi_2\!+\partial_{\overline{a}}\Phi_2\right) \cdot j \\
&= \left\{\frac{\overline{a}+a}{2\vert p\vert^{2}}-\frac{\left(a-\overline{a}\right) \left[2\vert p\vert^{2}-\left(a+\overline{a}\right)\overline{a}\right]}{2\vert p\vert^{2}\!\left[4 \vert p\vert^{2}-\left(a+\overline{a}\right)^{2}\right]}-\frac{\left(a-\overline{a}\right)\left[2\vert p\vert^{2}-\left(a+\overline{a}\right)a\right]}{2\vert p\vert^{2}\!\left[4 \vert p\vert^{2}-\left(a+\overline{a}\right)^{2}\right]}\right\} \\
& -\left\{\frac{b\left[2\vert p\vert^{2}-\left(a+\overline{a}\right)\overline{a}\right]}{\vert p\vert^{2}\left[4 \vert p\vert^{2}-\left(a+\overline{a}\right)^{2}\right]}+\frac{b\left[2\vert p\vert^{2}-\left(a+\overline{a}\right)a\right]}{\vert p\vert^{2}\left[4 \vert p\vert^{2}-\left(a+\overline{a}\right)^{2}\right]}\right\}\cdot j \\
&=\frac{a+\overline{a}}{2\vert p\vert^{2}}-\frac{\left(a-\overline{a}\right)\left[4\vert p\vert^{2}- \left(a+\overline{a}\right)^{2}\right]}{2\vert p\vert^{2}\left[4 \vert p\vert^{2}-\left(a+\overline{a}\right)^{2}\right]}-\frac{b\left[4\vert p\vert^{2}- \left(a+\overline{a}\right)^{2}\right]}{\vert p\vert^{2}\left[4 \vert p\vert^{2}-\left(a+\overline{a}\right)^{2}\right]}\cdot j \\&=\frac{\overline{a}}{\vert p\vert^{2}}- \frac{b}{\vert p\vert^{2}}\cdot j=\frac{\overline{p}}{\vert p\vert^{2}}=\frac{1}{p}.
\end{align*}

We see that the first derivative of the quaternionic function $\psi_{H}\!\left(p\right)=\text{Ln}\left(p\right)$ has the same form as the first derivative of the logarithmic function in complex (and real) analysis: $\left(\text{Ln}\left(\xi\right)\right)^{{\left(1\right)}}\!=~\!\frac{1}{\xi}.$ Since the function $\frac{1}{p}$ considered in Example  \ref{Exa3} above is H-holomorphic, the first derivative of the function
$\text{Ln}\left(p\right)$ is H-holomorphic too. The computation of the higher derivatives of $\text{Ln}\left(p\right)$ repeats the results of Example \ref{Exa3}. Thus the quaternion logarithmic function $\psi_{H}\!\left(p\right)= \text{Ln}\left(p\right)$ possesses H-holomorphic derivatives of all orders everywhere in $\mathbb{H}\setminus \!\left\{0\right\}\!.$
\subsection{Example of  commutative behaviour multiplication of H-holomorphic functions}  \label{Exa5} 
 \textit{The quaternionic multiplication of the quaternionic holomorphic functions} $\psi_{H}\!\left(p\right)=p^{2}$ \textit{and} $\psi_{H}\!\left(p\right)=e^{p}$  \textit{behaves as commutative.} Let the functions $f=f_{1}+f_{2}\cdot j$ and $g=g_{1}+g_{2}\cdot j$ be arbitrary quaternionic functions in the Cayley–Dickson construction. The quaternionic multiplication \cite{ks:hn} of these functions can be writen as follows:
 \begin{equation} \label{Eq8}
f\cdot g=\left(f_{1}+f_{2}\cdot j\right)\cdot\left(g_{1}+g_{2}\cdot j\right) =Re\left(f\cdot g\right)+Im\left(f\cdot g\right)\cdot j,
\end{equation}
where $$Re\left(f\cdot g\right)=f_{1}\,g_{1}-f_{2}\,\overline{g_{2}},\,\,\,\,\, Im\left(f\cdot g\right)=f_{2}\,\overline{g_{1}}+f_{1}\,g_{2}$$\\
are the designations of the parts of the quaternionic product. \\
\textit{The quaternionic product} $p^{2}\cdot e^{p}.$ In this case we have (see Examples  \ref{Exa1},  \ref{Exa2} ) $f_{1}=a^{2}-b\overline{b},\,\, f_{2}=\left(a+\overline{a}\right)b,\,\, g_{1}\!=2\beta\cos v+\frac{\beta\left(a-\overline{a}\right)\sin v}{v}$ and $g_{2}=\frac{2\beta\, b\sin v}{v}.$ Subsituting these into (\ref{Eq8}) we obtain the following expressions: $$Re\left(p^{2}\cdot e^{p}\right)=f_{1}\,g_{1}-f_{2}\,\overline{g_{2}}=\left(a^{2}-b\overline{b}\right)\left[2\beta\cos v+\frac{\beta\left(a-\overline{a}\right)\sin v}{v}\right]-\left(a+\overline{a}\right)b\frac{2\beta\, \overline{b}\sin v}{v},$$ $$Im\left(p^{2}\cdot e^{p}\right)=f_{2}\,\overline{g_{1}}+f_{1}\,g_{2}=\left(a+\overline{a}\right)b\left[2\beta\cos v+\frac{\beta\left(\overline{a}-a\right)\sin v}{v}\right]+\left(a^{2}-b\overline{b}\right)\frac{2\beta\, b\sin v}{v},$$
where we take into account (in the sequel also) that in accordance with (\ref{Eq5}) we have $\beta=\overline{\beta},\,\,v=\overline{v}.$ \\
\textit{The quaternionic product} $e^{p}\cdot p^{2}.$ In this case we have $f_{1}\!=2\beta\cos v+\frac{\beta\left(a-\overline{a}\right)\sin v}{v},\,\,f_{2}\!=\frac{2\beta b\sin v}{v},\,\, g_{1}=a^{2}-b\overline{b}$  and $g_{2}=\left(a+\overline{a}\right)b.$ Substituting these into (\ref{Eq8}) we obtain as follows: 
\begin{align*} %Muster c & ohne Number
Re\left( e^{p}\cdot p^{2}\right)=f_{1}\,g_{1}-f_{2}\,\overline{g_{2}}=\left[2\beta\cos v+\frac{\beta\left(a-\overline{a}\right)\sin v}{v}\right]\left(a^{2}-b\overline{b}\right)-\frac{2\beta\, b\sin v}{v}\left(a+\overline{a}\right)\overline{b}, \\ Im\left( e^{p}\cdot p^{2}\right)=f_{2}\, \overline{g_{1}}+f_{1}\,g_{2}= \frac{2\beta\, b\sin v}{v}\left(\overline{a}^{2}-b\overline{b} \right) +\left[2\beta\cos v+\frac{\beta\left(a-\overline{a}\right)\sin v}{v}\right]\left(a+\overline{a}\right)b.
\end{align*}
We see immediately that$$Re\left(p^{2}\cdot e^{p}\right)=Re\left( e^{p}\cdot p^{2}\right).$$

Now we verify that $Im\left(p^{2}\cdot e^{p}\right)=Im\left( e^{p}\cdot p^{2}\right),$ i.e.
\begin{align*} %Muster c & ohne Number
&\left(a+\overline{a}\right)b\left[2\beta\cos v+\frac{\beta\left(\overline{a}-a\right)\sin v}{v}\right]+\left(a^{2}-b\overline{b}\right)\frac{2\beta\, b\sin v}{v}\\  &= \frac{2\beta\, b\sin v}{v}\left(\overline{a}^{2}-b\overline{b} \right) +\left[2\beta\cos v+\frac{\beta\left(a-\overline{a}\right)\sin v}{v}\right]\left(a+\overline{a}\right)b.
\end{align*}
Opening the brackets and simplifying the expressions, we finally obtain the following identity: $$\left(a^{2}+\overline{a}^{2}\right)\!\beta\sin v=\left(a^{2}+\overline{a}^{2}\right)\!\beta\sin v,$$ that is, 
$$Im\left(p^{2}\cdot e^{p}\right)=Im\left( e^{p}\cdot p^{2}\right).$$
Thus, it is shown that the quaternionic multiplication of the quaternionic holomorphic functions $\psi_{H}\!\left(p\right)=p^{2}$ and $\psi_{H}\!\left(p\right)=e^{p}$ behaves as commutative.
 \subsection{Example of left and right quotients equality of H-holomorphic functions}  \label{Exa6}
 \textit{The left quotient of} $e^{p}$ \textit{by} $p$ \textit{is equal to the right one.} In accordance with \cite[p.20]{pm:eac} the left $X_{L}$ and and the right $X_{R}$ quotients in this case are correspondingly the following~: $X_{L}=p^{{-1}}\! \cdot e^{p}$ and $X_{R}=e^{p}\cdot p^{{-1}},\,\,\,\,p\neq 0.$ Given that the functions $p^{{-1}}$ (see Example of $p^{{-1}}$  \ref{Exa3}) and $e^{p}$ (see Example of $e^{p}$  \ref{Exa2}) are H-holomorphic and their quaternionic multiplication behaves as commutative, we could immediately conclude that $X_{L}=X_{R}$, however we will show this by direct computations. Using the above expressions for $\Phi_{1}$ and $\Phi_{2}$ of the functions $p^{{-1}}$ and $e^{p},$ and calculating quite similarly to the previous example  \ref{Exa5} we get, omitting intermediate results, the following expressions:
 \begin{align*} %Muster c & ohne Number
&Re\left(p^{{-1}}\!\!\cdot e^{p}\right)=\frac{\overline{a}}{a\overline{a}+b\overline{b}}\left[2\beta\cos v+\frac{\beta\left(a-\overline{a}\right)\sin v}{v}\right]+\frac{b}{a\overline{a}+b\overline{b}}\,\frac{2\beta\, \overline{b}\sin v}{v}, \\ &Im\left(p^{{-1}}\!\!\cdot e^{p}\right)=-\frac{b}{a\overline{a}+b\overline{b}}\left[2\beta\cos v+\frac{\beta\left(\overline{a}-a\right)\sin v}{v}\right]+\frac{\overline{a}}{a\overline{a}+b\overline{b}}\frac{2\beta\, b\sin v}{v}, \\ &Re \left(e^{p}\!\!\cdot p^{{-1}} \right)=\left[2\beta\cos v+\frac{\beta\left(a-\overline{a}\right)\sin v}{v}\right]\frac{\overline{a}}{a\overline{a}+b\overline{b}}+\frac{2\beta\,b\sin v}{v}\frac{\overline{b}}{a\overline{a}+b\overline{b}}\,,\\  &Im \left(e^{p}\!\!\cdot p^{{-1}} \right)=\frac{2\beta\,b\sin v}{v}\frac{a}{a\overline{a}+b\overline{b}}-\left[2\beta\cos v+\frac{\beta\left(a-\overline{a}\right)\sin v}{v}\right]\frac{b}{a\overline{a}+b\overline{b}}\,.
\end{align*}

We see immediately that$$ Re \left(p^{{-1}} \!\!\cdot e^{p} \right)=Re \left(e^{p}\!\!\cdot p^{{-1}} \right).$$

Now we show that $$ Im \left(p^{{-1}} \!\!\cdot e^{p} \right)=Im \left(e^{p}\!\!\cdot p^{{-1}} \right)\!,$$\\ that is, 
 \begin{align*} 
&-\frac{b}{a\overline{a}+b\overline{b}}\left[2\beta\cos v+\frac{\beta\left(\overline{a}-a\right)\sin v}{v}\right]+   \frac{\overline{a}}{a\overline{a}+b\overline{b}}\frac{2\beta\,b\sin v}{v} \\ &= \frac{2\beta\,b\sin v}{v}\frac{a}{a\overline{a}+b\overline{b}}- \left[2\beta\cos v+\frac{\beta\left(a-\overline{a}\right)\sin v}{v}\right]\frac{b}{a\overline{a}+b\overline{b}}\,.
\end{align*}
Opening the brackets and simplifying the expressions, we finally obtain the following identity:
\begin{align*} %Muster c & ohne Number
\overline{a}\,\,\frac{\beta\,b\sin v}{v}+\frac{\beta\,ba\sin v}{v}=\frac{\beta\,b\sin v}{v}a+\frac{\beta\,b\overline{a}\sin v}{v},
\end{align*}
that is, $$ Im \left(p^{{-1}} \!\!\cdot e^{p} \right)=Im \left(e^{p}\!\!\cdot p^{{-1}} \right).$$

We have shown that the left quotient of  $e^{p}$ \textit{by} $p$ is equal to the right one and at the same time confirmed once again that the quaternionic multiplication of H-holomorphic functions behaves as commutative.

The commutativity property of quternionic multiplication of H-holomorphic functions as well as the equality of the left and the right quotients of H-holomorphic functions can be illustrated with special computer programs \cite{pm:ph} when using various examples of H-holomorphic functions

As shown above, quaternionic generalization (\ref{Eq3}) of complex Cauchy-Riemann's equations is based on the requirement of uniqueness of a quaternionic derivative, which is needed to explore steady state vector fields in 3D space. It is not superfluous to note that the physical formulation of a problem played initially an important role in the theory of complex-differentiable functions, and the so-called complex Cauchy-Riemann equations were found \cite{re:fth} as early as in 1752 in d'Alembert's doctrine about planar fluid flow.
 
  \label{sec:intro}
\bibliography{bibliography}
\end{document}